\documentclass[11pt,a4paper]{article}
\usepackage{amssymb,amsmath}
\usepackage{graphicx}
\usepackage{url}
\usepackage{a4wide}
\usepackage{color}
\newcommand {\eq} [1] {\begin{equation}\label{#1}}
\newcommand {\en} {\end{equation}}
\newcommand {\eqn}      {\begin{eqnarray}}
\newcommand {\enn}      {\end{eqnarray}}
\newcommand {\bstar}    {\begin{eqnarray*}}
\newcommand {\estar}    {\end{eqnarray*}}
\newcommand {\B}        {{\bf B}}

\newcommand {\Q}        {{\bf Q}}

\newcommand {\T}        {{\bf T}}

\newcommand {\proof} {\par{\it Proof}. \ignorespaces}
\newcommand {\eproof}
      {\space
        {\ \vbox{\hrule\hbox{\vrule height1.3ex\hskip0.8ex\vrule}\hrule}}
        \par}

\def\u{u}

\newcommand {\corank}     {\mathop{\rm corank}\nolimits}

\newcommand {\mat}      [1] {\left[\begin{array}{#1}}
\newcommand {\rix}          {\end{array}\right]}
\newtheorem{theorem}            {Theorem}
\newtheorem{lemma}     [theorem]{Lemma}
\newtheorem{definition}[theorem]{Definition}

\newtheorem{example} [theorem]          {Example}
\newtheorem{remark} [theorem]           {Remark}

\newtheorem{hypothesis}  [theorem]        {Hypothesis}
\newcommand {\rank}       {\mathop{\rm rank}\nolimits}
\newcommand {\diag}    {\mathop{\rm diag}\nolimits}

\title{Port-Hamiltonian descriptor systems}

\author{Christopher Beattie\footnotemark[1] \, and Volker Mehrmann\footnotemark[2] \, and Hongguo Xu\footnotemark[3]\, and  Hans Zwart\footnotemark[4]
}
\begin{document}
\maketitle

\begin{abstract}
The modeling framework of port-Hamiltonian systems is systematically extended to constrained dynamical systems (descriptor systems, differential-algebraic equations). A new algebraically and geometrically defined system structure is derived. It is shown that this structure is invariant under equivalence transformations, and that it is adequate also for the modeling of high-index descriptor systems. The regularization procedure for descriptor systems to make them suitable for simulation and control is modified to deal with the port-Hamiltonian structure. The relevance of the new structure is demonstrated with several examples.
\end{abstract}


\maketitle
\renewcommand{\thefootnote}{\fnsymbol{footnote}}

\footnotetext[1]{
Department of Mathematics, Virginia Tech, Blacksburg, VA 24061, USA.
\texttt{beattie@vt.edu}. Supported by {\it Einstein Foundation Berlin},
through an Einstein Visiting Fellowship.}
\renewcommand{\thefootnote}{\arabic{footnote}}

\footnotetext[2]{
Institut f\"ur Mathematik MA 4-5, TU Berlin, Str. des 17. Juni 136,
D-10623 Berlin, FRG.
\texttt{mehrmann@math.tu-berlin.de}. Supported by {\it Einstein Foundation Berlin} via the Einstein Center ECMath and by Deutsche Forschungsgemeinschaft via Project A02 within CRC 1029 'TurbIn'}
\footnotetext[3]{
Department of Mathematics, University of Kansas, Lawrence, KS 66045, USA.
\texttt{xu@math.ku.edu}.
Partially supported by {\it Alexander von Humboldt Foundation}
and by {\it Deutsche Forschungsgemeinschaft},
through the DFG Research Center {\sc Matheon}
{\it Mathematics for Key Technologies} in Berlin.
}
\footnotetext [4] {University of Twente, Department of Applied Mathematics, P.O. Box 217, 7500 AE Enschede, The Netherlands.
\texttt{h.j.zwart@utwente.nl} and Eindhoven University of Technology, Department of Mechanical Engineering, P.O. Box 513,
5600 MB Eindhoven, The Netherlands.}


\section{Introduction}

Modeling packages such as  {\sc modelica} (\texttt{https://www.modelica.org/}), {\sc Matlab/Simu\-link} ({\texttt{http://www.mathworks.com}}) or {\sc Simpack} \cite{simpack} have come to provide excellent capabilities for the automated generation of models describing dynamical systems originating in different physical domains that may include mechanical, mechatronic, fluidic, thermic, hydraulic, pneumatic, elastic, plastic, or electric components  \cite{BalHPS05,EicF98,HilH06,Sch93,SchK01}.
Due to the explicit incorporation of constraints, the resulting systems comprise {\em differential-algebraic equations (DAEs)}, also referred to as  {\em descriptor systems} in the system theory context. Descriptor systems may contain hidden constraints, consistency requirements for initial conditions, and unexpected regularity requirements. Therefore, these models usually require further regularization to be suitable for numerical simulation and control, see \cite{CamKM11,KunM06,LamMT13}.
Our main focus will be on linear-time varying descriptor systems, as they may arise from the linearization of nonlinear DAE systems along a (non-stationary) reference trajectory, see \cite{Cam95}. These have the form
\begin{eqnarray}
    E(t) \dot x(t) & =&A(t) x(t)+B(t) u(t), \nonumber \\ 
 y &=& C(t) x(t)+D(t) u(t), \label{eq:lindesout}
\end{eqnarray}
together with an initial condition $x(t_0)=x_0$.  The coefficient matrices  $E,A \in C^0({\mathbb I},{\mathbb R}^{n,n})$, $B\in C^0({\mathbb I},{\mathbb R}^{n,m})$, $C\in C^0({\mathbb I},{\mathbb R}^{m,n})$, and $D\in C^0({\mathbb I},{\mathbb R}^{m,m})$, where we denote by $C^j({\mathbb I},\mathcal X)$ $j\in \{0,1,2,3,\ldots\}$ the set of $j$-times continuously differentiable functions from a compact time interval $ \mathbb I=[t_0,t_f]\subseteq \mathbb R$ to $\mathcal X={\mathbb R}^{n}$. If it is otherwise clear from the context, the argument~$t$ of the coefficient functions is suppressed.

An important development in recent years has been to employ \emph{energy based modeling} via bond graphs \cite{Bre08,CouJMTB08}.  This has been implemented recently in {\sc 20-sim} (\texttt{http://www.20sim.com/}), for example.   The resulting systems have a \emph{port-Hamiltonian  (pH) structure}, see e.~g. \cite{GolSBM03,JacZ12,OrtSMM01,Sch04,Sch06}, that encodes underlying physical principles such as conservation laws directly into the structure of the system model.  The standard form for {\em pH systems} appears as
 \begin{align}
\dot x&=\left(J-R\right)\nabla_{\!x}{\mathcal H}(x)+(B-P)u,\nonumber \\
y&=\ (B+P)^T \nabla_{\!x}{\mathcal H}(x)  + (S+N) u, \label{PHdef}
\end{align}
where the function ${\mathcal H}(x)$ is the \emph{Hamiltonian} which describes  the distribution of internal energy among energy storage elements of the system, $J=-J^T \in \mathbb R^{n,n}$ is the \emph{structure matrix} describing energy flux among energy storage elements within the system; $R=R^T\in \mathbb R^{n,n}$ is the \emph{dissipation matrix} describing energy dissipation/loss in the system;  $B\pm P\in\mathbb{R}^{n,m}$ are \emph{port} matrices,  describing the manner in which energy enters and exits the system, and $S+N$, with $S=S^T\in \mathbb R^{m,m}$ and $N=-N^T \in \mathbb R^{m,m}$, describes the direct \emph{feed-through} from input to output.  It is necessary that
\begin{equation} \label{Kdef}
W=\left[\begin{array}{lc}
R & P \\[1mm]
P^T & S
\end{array}\right] \geq 0,
\end{equation}
where we write $W>0$ (or $W\geq 0$) to assert that a real symmetric matrix $W$ is positive definite (or positive semi-definite).
Port-Hamiltonian systems generalize \emph{Hamiltonian systems}, in the sense that the \emph{conservation of energy} for Hamiltonian systems is replaced by the {\em dissipation inequality}:
\begin{equation}  \label{DissipIneq}
 {\mathcal H}(x(t_1))-{\mathcal H}(x(t_0)) \leq \int_{t_{0}}^{t_{1}} y(t)^Tu(t)\ dt.
\end{equation}
%
In the language of system theory, \eqref{DissipIneq} shows that the dynamical system described in (\ref{PHdef}) is a {\em passive} system \cite{ByrIW91}. Furthermore, ${\mathcal H}(x)$ defines a Lyapunov function for the unforced system, so pH systems are implicitly Lyapunov stable \cite{HinP05}.
Inequality (\ref{DissipIneq}) is an immediate consequence of (\ref{Kdef}) and holds even when the coefficient matrices $J$, $R$, $B$, $P$, $S$, and $N$ depend on $x$ or explicitly on time $t$, see \cite{MasSB92}, or when they are defined as linear operators acting on infinite dimensional spaces \cite{JacZ12,SchM02}.

The physical properties of pH systems are encoded in the algebraic structure of the coefficient matrices and in geometric structures associated with the flow of the differential equation. This leads to a remarkably robust modeling paradigm that greatly facilitates the combination and manipulation of pH systems.  Note in particular that the family of pH systems is closed under \emph{power-conserving interconnection} (see \cite{Kle13});  model reduction of pH systems via Galerkin projection yields (smaller) pH systems \cite{BeaG11,GugPBS12,PolS10}; and conversely, pH systems are easily extendable in the sense that new state variables can be included while preserving the structure of (\ref{PHdef}), and so, the range of application of the model can be increased while ensuring that the basic conservation principle \eqref{DissipIneq} remains in force.


When state constraints are included in a pH system, the resulting system is  a {\em port-Hamiltonian descriptor system (differential-algebraic equation)} (pHDAE). pHDAE systems arise also in singularly perturbed pH systems when small parameters are set to zero, see \cite{Sch13}.
Significantly, there is no systematic way that has yet emerged to describe this problem class consistently, in a way that reflects both the pH structure and the DAE structure accurately.  The first main topic of this paper is to propose such a systematic approach.  This is a challenging task, in particular when  constraints of the DAE are 'hidden', which is often signaled with the terminology 'high-index DAE' \cite{BreCP96,KunM06,LamMT13}.  Such DAEs are not well-suited for numerical simulation and control and so, either a reformulation or a regularization of the model must first be carried out, \cite{CamKM11,KunM06}. We will briefly summarize the fundamentals of this technique in Section~\ref{sec:nonl}.

It is sometimes stated in the literature, see e.~g. \cite{Sch13}, that port-Hamiltonian DAEs are of \emph{differentiation-index at most one}, i.~e., that they do not contain hidden constraints arising from derivatives. In contrast, we will show that higher-index pHDAEs  are actually very common and so a regularization procedure is necessary.  Unfortunately, the usual regularization strategies  do not preserve a given pHDAE structure of the model and so, how one should go about this task while respecting pHDAE structure is the second main topic of the paper.

The paper is organized as follows. In Section~\ref{phdaedef} we give a definition of \emph{port-Hamiltonian differential-algebraic systems} and demonstrate that this is a relevant class for many applications. The main properties of this new class of pHDAE systems (such as stability and dissipativity) are discussed in Section~\ref{sec:phprop}. Section~\ref{sec:nonl} extends the definition to the nonlinear case. The analysis of `index at most one' pHDAEs is discussed in Section~\ref{sec:index_one} while the structured regularization procedure is discussed in Section~\ref{sec:index_reduction}.

\section{Linear Port-Hamiltonian Differential-Algebraic Equations}\label{phdaedef}
In this section we introduce a new definition of systems of  port-Hamiltonian descriptor systems (pHDAEs). Our new definition is slightly different from the concepts discussed in \cite{Sch13} and is based on the concept of skew-adjoint differential-algebraic operators, see \cite{KunMS14} for the corresponding self-adjoint case.
\begin{definition} \label{def:skew-adjoint}
A (differential-algebraic) operator
\[
\mathcal L:= \mathcal E \frac{d}{dt}-\mathcal A : \Omega\subset C^1({\mathbb I},{\mathbb R}^{n})\to  C^0({\mathbb I},{\mathbb R}^{n})
\]
with coefficient functions $\mathcal E\in C^1(\mathbb I,\mathbb R^{n,n})$, $\mathcal A\in C(\mathbb I,\mathbb R^{n,n})$
is called {\em skew-adjoint},
if $\mathcal E^T(t)=\mathcal E(t)$ and $\dot {\mathcal E}(t)=-(\mathcal A(t)+\mathcal A^T(t))$ for all $t\in \mathbb I$.
\end{definition}
This definition is motivated by the following observation: starting with vector functions $x_1(t),x_2(t)$ that are absolutely continuous on the interval
${\mathbb I}=(t_0, t_f)$ each with square integrable derivative and $x_i(t_0)=x_i(t_f)=0$ for $i=1,2$, and then denoting the usual $L_2$ inner product
as $\langle x_1,x_2 \rangle=\int_{t_0}^{t_f} x_2^Tx_1\, dt$, we have
\begin{eqnarray*}
\langle x_1, \mathcal L (x_2) \rangle &=& \langle x_1, \mathcal  E \dot x_2-\mathcal A x_2\rangle
=\langle  x_1, \frac{d}{dt} (\mathcal E x_2)- \mathcal A x_2 -\dot {\mathcal E} x_2 \rangle\\
&=&  x_2^T\mathcal E x_1|_{t_0}^{t_f}-\langle  \mathcal E^T \dot x_1,x_2\rangle -\langle (\mathcal A^T+\dot {\mathcal E}^T)x_1,x_2\rangle\\
&=& \langle  -\mathcal E^T \dot x_1 -(\mathcal A^T+\dot {\mathcal E}^T) x_1,x_2\rangle =  \langle  -\mathcal E \dot x_1 +\mathcal A x_1,x_2\rangle.
\end{eqnarray*}
%
%
%

So formally, the adjoint operator $\mathcal L^*$ satisfies $\mathcal L^*=-\mathcal L$.
Note the boundary terms arising in partial integration will vanish under a wide variety of conditions replacing the requirement of zero end conditions on $x_1(t)$ and $x_2(t)$.
\begin{remark}\label{rem:skew-sym}{\rm
In the context of densely defined unbounded operators, recall that \emph{symmetric} operators are those with adjoints that are extensions of the original operator, so one might use analogously the terminology \emph{skew-symmetric} operator instead of skew-adjoint operator here. To be consistent with the terminology in \cite{KunMS14} where  self-adjoint DAE operators were introduced, we prefer to use its natural cousin, \emph{skew-adjoint} operator.}
\end{remark}
Skew-adjoint operators stay skew-adjoint under time-varying congruence transformations.
\begin{lemma}\label{lem:skewa}
Consider a skew-adjoint differential-algebraic operator
\[
\mathcal L:= \mathcal E \frac{d}{dt}-\mathcal A : \Omega\subset C^1({\mathbb I},{\mathbb R}^{n})\to  C^0({\mathbb I},{\mathbb R}^{n})
\]
with
coefficient functions $\mathcal E\in C^1(\mathbb I,\mathbb R^{n,n})$ and $\mathcal A\in C(\mathbb I,\mathbb R^{n,n})$.
Then for every $\mathcal V\in C^1(\mathbb I,\mathbb R^{n,r})$,  the operator
$\mathcal L_{\mathcal V}$ defined by
\[
\mathcal L_{\mathcal V}(x):= {\mathcal V}^T \mathcal E \mathcal V \dot x -
(\mathcal V^T \mathcal A \mathcal V -\mathcal V^T \mathcal E \dot {\mathcal V})x
\]
is again skew-adjoint, i.~e., $\mathcal L_{\mathcal V}^*=-\mathcal L_{\mathcal V}$.
\end{lemma}
\proof
Since ${\mathcal V}^T \mathcal E \mathcal V=({\mathcal V}^T \mathcal E \mathcal V)^T $,
it remains to consider the coefficient of $x$.
Using $\mathcal E^T=\mathcal E$ and $\dot{\mathcal E}=-(\mathcal A+\mathcal A^T)$,
we have
\begin{eqnarray*}
\frac{d}{dt}(\mathcal V^T\mathcal E\mathcal V)&=&
\dot{\mathcal V}^T\mathcal E\mathcal V+
\mathcal V^T\dot{\mathcal E}\mathcal V+\mathcal V^T\mathcal E\dot{\mathcal V}\\
&=&\dot{\mathcal V}^T\mathcal E\mathcal V-
\mathcal V^T({\mathcal A}+\mathcal A^T)\mathcal V+\mathcal V^T\mathcal E\dot{\mathcal V}
\\
 &=&
-(\mathcal V^T \mathcal A \mathcal V-\mathcal V^T\mathcal E\dot{\mathcal V})
-(\mathcal V^T \mathcal A \mathcal V-\mathcal V^T\mathcal E\dot{\mathcal V})^T.
\end{eqnarray*}
\eproof
It should be noted that for any $t\in \mathbb I$ and $x\in C(\mathbb I,\mathbb R^{n})$ we have ${\mathcal L}_{\mathcal V}(x(t))={\mathcal V}^T(t){\mathcal L}(\mathcal V(t)x(t))$.
\begin{remark}\label{rem:proj}{\rm
Note that in Lemma~\ref{lem:skewa}, ${\mathcal V }$ need be neither invertible nor square, and in particular a time-varying compression ${\mathcal V }=\mat{c} I_r \\ P(t) \rix $ will produce a permissable skew-adjoint operator.
}
\end{remark}
\begin{definition}\label{def:pHDAE}
A linear variable coefficient descriptor system of the form
\begin{eqnarray}
E \dot x &=&\left [(J-R) Q -E K\right ]x + (B-P)u, \nonumber\\
y&=& (B+P)^T Q x  + (S+N) u, \label{pHDAE}
\end{eqnarray}
with $E,Q\in C^1({\mathbb I},{\mathbb R}^{n,n})$, $J,R, K \in C^0({\mathbb I},{\mathbb R}^{n,n})$,
$B,P\in C^0({\mathbb I},{\mathbb R}^{n,m})$,  $S=S^T, N=-N^T\in C^0({\mathbb I},{\mathbb R}^{m,m})$
is called  \emph{port-Hamiltonian descriptor system (port-Hamiltonian differential-algebraic system) (pHDAE)} if the following properties are satisfied:
\begin{itemize}
\item [i)] the differential-algebraic operator
\begin{equation}
\label{skewop}
\mathcal L:= Q^TE \frac{d}{dt} -(Q^TJ Q -Q^TE K):  {\mathcal D}\subset C^1({\mathbb I},{\mathbb R}^{n})\to  C^0({\mathbb I},{\mathbb R}^{n})
\end{equation}
is \emph{skew-adjoint}, i.~e. we have that  $Q^TE \in C^{1}(\mathbb I, \mathbb R^{n,n})$ and for all $t\in \mathbb I$,
\begin{eqnarray*}
Q^T(t)E(t)&=&E^T(t)Q(t), \mbox{ and  }\\
 \frac{d}{dt}( Q^T(t) E(t)) &=&Q^T(t)[E(t) K(t) - J(t)Q(t)] +[E(t) K(t) - J(t)Q(t)]^TQ(t);
\end{eqnarray*}
 \item[ii)] the \emph{Hamiltonian} function defined as
\begin{equation}\label{defham}
\mathcal H(x): = \frac 12 x^TQ^TEx:  C^1({\mathbb I},{\mathbb R}^{n})\to C^1({\mathbb I},\mathbb R)
\end{equation}
is bounded from below by a constant, i.~e., ${\mathcal H}(x(t))\geq h_0\in \mathbb R$ uniformly for all $t\in \mathbb I$ and all solutions $x$ of (\ref{pHDAE});
%
\item [iii)] the matrix function
\begin{equation} \label{Wdef}
W:=\left[\begin{array}{lc}
Q^T R Q& Q^T P \\[1mm]
P^T Q & S
\end{array}\right] \in C^0({\mathbb I},\mathbb{R}^{n+m,n+m})
\end{equation}
is positive semidefinite, i.~e., $W(t)=W^T(t)\geq 0$  for all $t\in \mathbb I$.
\end{itemize}
%
\end{definition}
Assumption ii) in Definition~\ref{def:pHDAE} can be formulated as a condition on the matrix function $Q^TE$ as shown in the following Lemma.
\begin{lemma}  Assumption ii) in Definition \ref{def:pHDAE} is equivalent to the assertion that $Q^T(t)E(t)$ is positive semidefinite for all $t\in \mathbb I$.
\end{lemma}
\proof
Evidently, if $Q^T(t)E(t)$ is positive semidefinite for all $t\in \mathbb I$, then ${\mathcal H}(x)$ is semibounded with $h_0=0$.
Now suppose that ${\mathcal H}(x)$ is semibounded, say with $h_0<0$, but $Q^T(\hat{t})E(\hat{t})$ fails to be positive semidefinite at some time point  $\hat{t}$.  Then there exists a state $\hat{x}$ such that ${\mathcal H}(\hat{x}(\hat{t}))<0$.  By scaling $\hat{x}$ by $\kappa>\sqrt{\frac{h_0}{   {\mathcal H} ( \hat{x}(\hat{t}))}}$, we find that ${\mathcal H}(\kappa\hat{x}(\hat{t}))=\kappa^2 {\mathcal H}(\hat{x}(\hat{t}))<h_0$ giving a contradiction, so $Q^T(t)E(t)$ must indeed be positive semidefinite for all $t\in \mathbb I$.
\eproof
Besides the matrix function $E$ in front of the derivative and the different definition of the Hamiltonian, which gives the option of having singular matrices $E$ and $Q$, a major difference to the definition of standard pH systems is the extra additive term $-EKx$ on the right hand side of (\ref{pHDAE}), which is needed to accommodate time-varying changes of basis.

Note further that in this definition no further properties of the differential-algebraic operator are assumed. In particular it is not assumed that it has a certain index as a differential-algebraic equation.
\begin{example}\label{ex:circuit}{ \rm
Consider the model of a simple RLC network, see e.~g. \cite{Dai89,Fre11}, given by a linear constant coefficient DAE
\begin{eqnarray}\label{eq:RLC_network_1}
\underbrace{\mat{ccc} G_c C G_c^T & 0 & 0\\ 0 & L & 0 \\ 0 & 0 & 0 \rix}_{:=E}
\mat{c}\dot{V} \\ \dot{I_l}\\ \dot{I_v} \rix=
\underbrace{\mat{ccc} -G_r R_r^{-1}G_r^T & -G_l & -G_v \\ G_l^T & 0 & 0\\ G_v^T & 0 & 0 \rix}_{:=(J-R)I}
\mat{c}V \\ I_l\\ I_v \rix,
\end{eqnarray}
with real symmetric constant matrices $L > 0$, $C > 0$, $R_r> 0$ describing inductances,  capacitances, and resistances, respectively that are present in the network.  Here, $G_v$ is of full column rank, and the subscripts
$r,\,c,\,l,$ and $v$  refer to edge quantities corresponding to the resistors,
capacitors, inductors, and voltage sources,
while $V$, $I$ denote the voltage and current, respectively, on or across the branches of the given RLC network. This model has a  pHDAE structure  with vanishing $B,P,S,N,K$, the matrix $Q$ is the identity,  $E=E^T$, $J=-J^T$, $Q^TRQ=R\geq 0$, and
\[\mathcal H= \mat{c}V \\ I_l \\ I_v \rix^T E \mat{c}V \\ I_l \\ I_v \rix=
\mat{c}V \\ I_l\rix^T\mat{ccc} G_c C G_c^T & 0 \\ 0 & L  \rix \mat{c}V \\ I_l  \rix.
\]
}
\end{example}

\begin{example}\label{ex:gas}
{\rm In \cite{EggK16,EggKLMM17} the  propagation of pressure waves on the acoustic time scale in a network of  gas pipelines is considered and an infinite-dimensional pHDAE is derived. A structure preserving mixed finite element space discretization leads to a block-structured pHDAE system
\begin{eqnarray}
E\dot x &=&(J-R)Qx+Bu, \nonumber\\ 
y&=&B^TQ x \label{eq:ph1}, \\
x(t_0)&=& x^0, \nonumber 
\end{eqnarray}
with $Q=I$, $P=0$, $S+N=0$,
\[
E=\begin{bmatrix} M_1 & 0 & 0 \\ 0 & M_2 & 0 \\ 0 & 0 & 0 \end{bmatrix},
J=\begin{bmatrix} 0 & -\tilde  G & 0 \\ \tilde  G^T & 0 & \tilde N^T \\ 0 & -\tilde N & 0\end{bmatrix}, R=\begin{bmatrix} 0 & 0 & 0 \\ 0 & \tilde  D & 0 \\ 0 & 0 & 0\end{bmatrix},
B=\begin{bmatrix} 0 \\ \tilde B_2 \\ 0 \end{bmatrix}, x=\begin{bmatrix}x_1 \\ x_2 \\ x_3 \end{bmatrix},
\]
where the vector valued functions $x_1: \mathbb R\to \mathbb R^{n_1},x_2: \mathbb R\to \mathbb R^{n_2}$ represent the discretized pressure and flux, respectively, and $x_3: \mathbb R\to \mathbb R^{n_3}$ represents the Lagrange multiplier for satisfying the space-discretized constraints. The coefficients  $M_1=M_1^T$, $M_2=M_2^T$, and $\tilde D=\tilde D^T$  are positive definite, and the matrices $\tilde N$ and $\begin{bmatrix} \tilde G^T & \tilde N^T \end{bmatrix}^T$ have full row rank. The Hamiltonian is given by $\mathcal H(x)=x^TE^TQx= x_1^T M_1 x_1+x^T_2 M_2 x_2$.}
\end{example}
Definition~\ref{def:pHDAE} brings the pH modeling framework and the DAE framework together in a structured way.
It should be noted, however, that in a DAE we may have hidden constraints that arise from differentiations, which are not explicitly formulated
and the formulation of the DAE that is used in simulation and control is not unique. One can for example add derivatives of constraints which leads to an over-determined system, then one can add dummy variables or Lagrange multipliers to make the number of variables equal to the number of equations or one can remove some of the dynamical equations to achieve this goal, see \cite{BreCP96,EicF98,KunM06,LamMT13} for detailed discussions on this topic.
To rewrite  these different formulations in the pHDAE formulation is not always obvious. Let us demonstrate this with an example from multi-body dynamics.
\begin{example}\label{ex:robot}{\rm A benchmark example for a nonlinear DAE system is the model of a two-dimensional three-link mobile manipulator, see \cite{BunBMN99,Hou94a}, which is modeled as
\begin{eqnarray}
     \tilde M(\Theta){\ddot \Theta}+\tilde C(\Theta,{\dot \Theta})+\tilde G(\Theta)&=&\tilde B_1\tilde u+\Psi^T\lambda, \nonumber\\
	\psi(\Theta) & = & 0, \label{eq:mbnon}
\end{eqnarray}
where $\Theta=\mat{ccc} \Theta_1& \Theta_2 & \Theta_3 \rix^T$
is the vector of joint displacements, $\tilde u$ is  vector of control
torques at the joints, $\tilde M$ is  mass
matrix, $\tilde C$ is the vector of centrifugal and Coriolis forces,
and $\tilde G$ is the gravity vector. The term  $\Psi^T\lambda$ with $\Psi = \frac{\partial \psi}{\partial \Theta}$ is the generalized constraint force with Lagrange multiplier $\lambda$  associated with the constraint
\[
\psi(\Theta)=\mat{c}
 l_1 \cos ({\Theta_1})+l_2 \cos(\Theta_1 + \Theta_2) +l_3
\cos(\Theta_1 + \Theta_2 + \Theta_3) l_3 -l\\  \Theta_1 + \Theta_2 +
\Theta_3
\rix=0.
\]
Besides the explicit constraint, this system contains the first and second time derivative of $\psi$ as hidden algebraic constraints, see e.~g. \cite{EicF98,KunM06}. There are several regularization procedures that one can employ to make the system better suited for numerical simulation and control. One possibility is to replace the original constraint by its time derivative $\Psi(\Theta) \dot \Theta=0$. In this case the model equation can easily be written in a pHDAE formulation. Using Cartesian coordinates for positions $p$,  scaling the constraint equation by $-1$, and linearizing around a stationary reference solution yields a linear constant coefficient DAE of the form
\begin{eqnarray}
	\tilde M \delta \ddot {p}  &=& - \tilde D\delta \dot {p} - \tilde S \delta {p}+\tilde G^T \delta \lambda+\tilde B_1
	\delta u,  \nonumber\\
	0 & = &  -\tilde G\delta \dot p ,\label{eq:mblin}
\end{eqnarray}
with pointwise symmetric positive definite matrices $\tilde M,\tilde S$  and pointwise symmetric and positive semidefinite $\tilde D$.
Adding a tracking output  of the form $y=\tilde B_1^T \delta \dot {p}$, see e.~g. \cite{HouM94}, and transforming to first order by introducing
\[
x=\mat{c} x_1 \\ x_2 \\ x_3  \rix:=\mat{c} \delta {\dot p} \\  \delta {p} \\ \delta \lambda  \rix,\ u=\delta  u,
\]
one obtains a linear constant coefficient pHDAE system
$E\dot x =(J-R)Q x+B u$, $y=B^T Q x$, with
\begin{eqnarray*}
E&:=&\mat{ccc} \tilde M & 0 & 0 \\ 0 & I & 0 \\ 0 & 0 & 0 \rix,\ R:=\mat{ccc} \tilde D &  0 & 0 \\
0 & 0 & 0 \\ 0 &0  & 0 \rix,\ Q:= \mat{ccc} I & 0 & 0 \\ 0 & \tilde S & 0 \\ 0 & 0 & I \rix,\\
J &:=&  \mat{ccc} 0 & -I & \tilde G^T \\ I & 0 & 0 \\ -\tilde G & 0 & 0 \rix,\ B:=\mat{c} \tilde B_1 \\ 0 \\ 0 \rix,\ P=0,\ S+N=0.
\end{eqnarray*}
The Hamiltonian in this case is given by $\mathcal H(x)= \mat{c} x_1 \\ x_2 \rix^T \mat{cc} \tilde M & 0  \\ 0 & \tilde S \rix\mat{c} x_1 \\ x_2 \rix$.

Since the Lagrange multipliers in the multibody framework can be interpreted as external forces, it is also possible to incorporate them in the input $(B-P)u$ to achieve a pHDAE formulation as in Definition~\ref{def:pHDAE}, but also other formulations are possible, e.~g. we can keep the original algebraic constraint as well and use an extra Lagrange multiplier for the first time derivative.
}
\end{example}

%
%
\begin{remark} \label{remark:pHDAEnoE}
A special case of \eqref{pHDAE} takes the following form:
\begin{eqnarray}
E \dot x &=&\ (J-R)x + (B-P)u, \nonumber\\
y&=& (B+P)^T x  + (S+N) u, \label{pHDAEnoE}
\end{eqnarray}
where $E=E^T\in C^1({\mathbb I},{\mathbb R}^{n,n})$, $A,R=R^T, K \in C^0({\mathbb I},{\mathbb R}^{n,n})$,
$B,P\in C^0({\mathbb I},{\mathbb R}^{n,m})$,  $S=S^T, N=-N^T\in C^0({\mathbb I},{\mathbb R}^{m,m})$ as before
but now we require,
\begin{itemize}
\item [i)] the differential algebraic operator
\begin{equation} \label{skewop-spec}
\mathcal L:= E \frac{d}{dt} -J :  D\subset C^1({\mathbb I},{\mathbb R}^{n})\to  C^0({\mathbb I},{\mathbb R}^{n})
\end{equation}
is \emph{skew-adjoint}, so that we have for all $t\in \mathbb I$,
\begin{eqnarray*}
 \frac{d}{dt} E(t) &=& - \left[J(t) + J(t)^T\right];
\end{eqnarray*}
\item [ii)] $E(t)$ is positive semidefinite: $E(t) \geq 0 $ for all $t\in \mathbb I$; and
\item [ii)]
$\displaystyle W(t):=\left[\begin{array}{lc}
 R(t)   &  P(t) \\[1mm]
P^T(t) &  S(t)
\end{array}\right] \geq 0$  for all $t\in \mathbb I$.
\end{itemize}
The effective \emph{Hamiltonian} is now
\begin{equation}\label{defhameff}
\mathcal H(x): = \frac 12 x^TE x:  C^1({\mathbb I},{\mathbb R}^{n})\to \mathbb R.
\end{equation}
Notice that in this model description
we have merged the roles of $Q$ and $E$. This is always possible when $Q$ is pointwise invertible, see Section~\ref{sec:phprop} but this formulation may not be possible when $Q$ is singular.
\end{remark}
%
%
\section{Properties of linear pHDAE systems}\label{sec:phprop}
To analyze the properties  of linear pHDAE systems, we first  derive the conservation of energy and the dissipation inequality.

\begin{theorem}\label{thm:dissin}
Consider the linear time-varying system (\ref{pHDAE}) and assume that this system satisfies condition {\em i)} of Definition \ref{def:pHDAE}. Then its (classical) solutions satisfy
\begin{equation}
\label{eq:bal}
  \frac{d}{dt} {\mathcal H}(x) =u^Ty -\mat{c}x\\ u \rix^T W \mat{c}x\\ u\rix,
\end{equation}
where $W$ is defined in (\ref{Wdef}).

Furthermore we have the  following properties.
\begin{itemize}
\item [i)] If $W\equiv 0$, then $\frac{d}{dt} \mathcal H =u^Ty$.
%
\item [ii)]  If $W\ge 0$ for all $t\in{\mathbb I}$,  then the system satisfies the dissipation inequality
(\ref{DissipIneq}).
\end{itemize}
\end{theorem}
\proof
By Definition~\ref{def:pHDAE} 
we have
%
\begin{eqnarray*}
\frac{d}{dt} \mathcal H &=&\frac 12\left[ \dot x^T (Q^TE)x +x^T\frac{d}{dt} (Q^TE)x +x^T (Q^TE)\dot x\right]\\
&=& \frac 12 x^T\frac{d}{dt} (Q^TE)x +x^T Q^T(E \dot x)\\
&=&\frac 12 x^T\frac{d}{dt} (Q^TE)x +x^TQ^T([JQ-RQ-EK]x+Bu-Pu)\\
&=&\frac 12 x^T\frac{d}{dt} (Q^TE)x +x^TQ^TJQx-x^TQ^TRQx-x^TQ^TEKx+x^TQ^TPu+u^TB^TQx\\
&=&\frac 12 x^T\frac{d}{dt} (Q^TE)x +x^TQ^TJQx-x^TQ^TRQx-x^TQ^TEKx-x^TQ^TPu\\
&& \qquad+u^T(y-P^TQx-Su-Nu)\\
&=&u^Ty+\frac 12 x^T\frac{d}{dt} (Q^TE)x +x^TQ^TJQx-x^TQ^TRQx-x^TQ^TEKx\\
&& \qquad -x^TQ^TPu
-u^TP^TQx-u^TSu\\
&=& u^Ty+\frac 12 \left(x^T\frac{d}{dt} (Q^TE)x +x^T[Q^T(JQ-EK)+(JQ-EK)^TQ]x\right)\\
&& \qquad -\mat{c}x\\ u\rix^T W\mat{c}x\\ u\rix.
\end{eqnarray*}
 From the skew-adjointness of $\mathcal L$ we then have that
\[
\frac{d}{dt} \mathcal H =u^Ty-\mat{c}x\\ u\rix^T W\mat{c}x\\ u\rix.
\]
Part i) then follows immediately from the assumption $W\equiv0$,
while in Part ii) the fact that $W(t)\geq 0$ for all $t\in \mathbb I$ implies
that for any $t_1\ge t_0$,
\[
{\mathcal H}(x(t_1))-{\mathcal H}(x(t_0))
=\frac 12\int_{t_0}^{t_1}{d\over dt}{\mathcal H}\, dt
\le \int_{t_0}^{t_1}y^Tu\,dt.
\]
\eproof
An important feature of pHDAE systems is that a change of basis and a scaling with an invertible matrix function preserves the pHDAE structure and the Hamiltonian. 
\begin{theorem}\label{pHDAEinv}
Consider a pHDAE system of the form (\ref{pHDAE})  with  Hamiltonian
(\ref{defham}). Let $U\in C^0({\mathbb I},{\mathbb R}^{n,n})$ and $V\in C^1({\mathbb I},{\mathbb R}^{n,n})$ be pointwise invertible in $\mathbb I$.
Then the transformed DAE
\begin{eqnarray*}
\tilde E  \dot {\tilde x}  &=&[(\tilde J -\tilde R )\tilde Q -\tilde E  \tilde K  ]
\tilde x  + (\tilde B -\tilde P )u, \\
y &=& (\tilde B+\tilde P)^T\tilde Q  \tilde x   + (S +N ) u,
\end{eqnarray*}
with
\begin{eqnarray*}
\tilde E &=&U^{T} E V , \quad
\tilde Q =U^{-1} Q V ,\quad
\tilde J =U^T J U ,\\
\tilde R &=&U^T R U ,\quad
\tilde B =U^T B ,\qquad
\tilde P =U^T P ,\\
\tilde K &=&V^{-1} K V  +{V}^{-1} \dot{V} ,\qquad\quad
x =V \tilde x
\end{eqnarray*}
is still a pHDAE system with the same Hamiltonian $\tilde{\mathcal H}(\tilde x)=\frac{1}{2}\tilde x^T  \tilde Q^T \tilde E \tilde x ={\mathcal H}(x)$.
\end{theorem}
\proof 
The transformed DAE system is obtained from the original DAE system by setting  $x=V\tilde x$ in
(\ref{pHDAE}), by pre-multiplying with $U^T$, and by inserting $U U^{-1}$ in front of $Q$.
The transformed operator corresponding to $\mathcal L$ in (\ref{skewop}) is
\[
\mathcal L_V:=\tilde Q^T\tilde E\frac{d}{dt}-\tilde Q^T(\tilde J\tilde Q-\tilde E\tilde K).
\]
Because
\[
\tilde Q^T\tilde E=V^TQ^TEV,\quad
\tilde Q^T\tilde J\tilde Q=V^TQ^TJQV,
\quad
\tilde Q^T\tilde EV^{-1}\dot V
=V^T Q^TE\dot V,
\]
by Lemma~\ref{lem:skewa}, $\mathcal L_V$ is skew-adjoint, since $\mathcal L$
defined in (\ref{skewop}) is skew-adjoint. Hence,
\begin{eqnarray*}
\tilde Q^T\tilde E&=&\tilde E^T\tilde Q,\\
\frac{d}{dt}(\tilde Q^T\tilde E)&=&-\tilde Q^T(\tilde J \tilde Q-\tilde E \tilde K)-
(\tilde J\tilde Q-\tilde E \tilde K)^T\tilde Q.
\end{eqnarray*}
It is straightforward to show that $\tilde{\mathcal H}(\tilde x)={\mathcal H}(x)$ and
\[
\frac{d}{dt}{\tilde{\mathcal H}}(\tilde x)=y^Tu-\mat{c}\tilde x\\ u\rix^T\tilde W\mat{c}\tilde x\\ u\rix,
\]
where
\begin{eqnarray*}
\tilde W&=&\mat{cc}\tilde Q^T\tilde R\tilde Q&\tilde Q^T\tilde P\\
\tilde P^T\tilde Q&S\rix
=\mat{cc}V^TQ^TRQV&V^TQ^TP\\P^TQV&S\rix\\
&=&\mat{cc}V&0\\0&I\rix^T W\mat{cc}V&0\\0&I\rix,
\end{eqnarray*}
and $W$ is defined in (\ref{Wdef}). Because $W(t)$ is positive semidefinite for all $t\in \mathbb I$, so is $\tilde W(t)$.
Therefore, for any $t_1\ge t_0$,
\[
\tilde{\mathcal H}(\tilde x(t_1))-\tilde{\mathcal H}(\tilde x(t_0))\le \int_{t_0}^{t_1}y^T(t)u(t)dt,
\]
which establishes the dissipation inequality.
%
%
%
%
%
An important point to note is that the Hamiltonian stays invariant under time-varying changes of basis and the operator $\mathcal L_V$, the Hamiltonian $\tilde {\mathcal H}(\tilde x)$, and the matrix function $\tilde W$ are independent of
the choice of the matrix function $U$.

\begin{lemma}\label{pHDAErev}
Consider a pHDAE system
\begin{eqnarray*}
\tilde E  \dot {\tilde x}  &=&[(\tilde J -\tilde R )\tilde Q -\tilde E \tilde K )]
\tilde x  + (\tilde B -\tilde P )u, \\
y &=& (\tilde B+\tilde P)^T\tilde Q  \tilde x     +(S +N ) u
\end{eqnarray*}
with Hamiltonian $\tilde{\mathcal H}(\tilde x)=\frac12\tilde x^T \tilde Q^T \tilde E \tilde x $,
where $\tilde K \in C({\mathbb I},{\mathbb R}^{n,n})$.
%
If  $V_{\tilde K} \in C^1({\mathbb I},{\mathbb R}^{n,n})$ is a pointwise invertible solution of the matrix differential equation
$\dot V =V \tilde K $ (e. g. with the initial condition $V(t_0)=I$), then defining
\begin{eqnarray*}
 E &=&\tilde E V_K^{-1} ,\quad Q =\tilde Q V_k^{-1}, \\
 J &=&\tilde J ,\qquad \,\,\,\,\,R =\tilde R ,\qquad B =\tilde B, \\
 P &=&\tilde P ,\qquad\,\,\,\,\, \tilde x =V_K^{-1} x ,
\end{eqnarray*}
the system
\begin{eqnarray*}
E  \dot { x}  &=&( J - R ) Q
x  + (B -P )u, \\
y &=& (B+ P)^TQ  x   + (S +N ) u
\end{eqnarray*}
is again pHDAE with the same Hamiltonian ${\mathcal H}(x)=\tilde{\mathcal H}(\tilde x)=
\frac12 x^T Q^T E x $.
\end{lemma}
\proof
For a given matrix function $\tilde K$, the system $\dot V=V\tilde K$
always has a solution $V_K$ that is pointwise invertible.
The remainder of the proof follows by reversing the proof of Theorem~\ref{pHDAEinv} with $U=I$ and using that
${\dot V}_K V_K^{-1}= -V_K\frac{d}{dt} (V_K^{-1})$.
\eproof
Note again that if $K$ is real and skew-symmetric, then the matrix function $V_K$ in Lemma~\ref{pHDAErev} can be chosen to be pointwise real orthogonal.

\begin{remark}\label{rem:id}{\rm
Following Theorem~\ref{pHDAEinv},
if $E$ is pointwise invertible, then the original system can be transformed into the one with new
$\hat E$ being the identity, so into a standard port-Hamiltonian system; and whenever $Q$ is pointwise invertible, then the original system can be transformed into the one with new $\hat Q$ being the identity, see Remark~\ref{remark:pHDAEnoE}.
Which of these formulations is preferable will depend on the sensitivity (conditioning) of these
transformations. In the context of numerical simulation and control methods, these transformations
should be avoided if they are  ill-conditioned.
}
\end{remark}

\section{Nonlinear DAEs and pHDAEs}\label{sec:nonl}
In this section we briefly recall the theory of  nonlinear DAE systems and then extend these results to pHDAEs. Consider a general descriptor system of the form
%
\begin{eqnarray}
  F ( t, x,\dot x, u) &=&0,\nonumber \\
  x(t_0)&=&x^0 \nonumber \\
  \label{nledestate}
 y&=& G(t,x,u).  
\end{eqnarray}
%
Assume that
${F}\in C^0({\mathbb I}\times{\mathbb D}_x\times{\mathbb D}_{\dot x}\times{\mathbb D}_u,{\mathbb R}^n)$
and ${G}\in C^0({\mathbb I}\times{\mathbb D}_x\times{\mathbb D}_u,{\mathbb R}^m)$ are sufficiently smooth, and that
${\mathbb D}_x,{\mathbb D}_{\dot x}\subseteq{\mathbb R}^n$, and
${\mathbb D}_{u}$
are open sets.
Note that (in order to deal with pHDAEs) in contrast to the more general case in \cite{CamKM11}, we assume square systems with an equal number of equations and variables and with an equal number of inputs and outputs.

For the analysis and the regularization procedure we make use of the \emph{behavior approach} \cite{PolW98}, which introduces a  \emph{descriptor vector} $v=[x^T,u^T]^T$. We could also include the output vector $y$ in $v$, but in the context of pHDAEs it is preferable to keep the output equation separate.
The behavior formulation has the form
\begin{equation}\label{nlcontrol}
{\mathcal F} ( t, v,\dot v)=0,
\end{equation}
with  ${\mathcal F}\in C^0({\mathbb I}\times{\mathbb D}_v\times{\mathbb D}_{\dot v},{\mathbb R}^{n})$  together with a set of initial conditions $c(v(t_0))=v^0$ which results from the original initial condition. Note that although no initial condition is given for $u$ in the context of the regularization procedure discussed in \cite{CamKM11} such conditions may arise, so we formally state a condition for $v(t_0)$.

To  regularize DAEs for numerical simulation and control, see \cite{Cam87a,CamKM11,KunM06}, one uses the behavior system (\ref{nlcontrol})
and some or all of its derivatives to produce an equivalent  system with the same solution set (all variables keep their physical interpretation), but where all explicit and hidden constraints are available. The approach of \cite{CamKM11} (adapted for the analysis of pHDAEs) uses the state equation of (\ref{nlcontrol}) to form a derivative array, see  \cite{Cam87a},
\begin{equation}\label{infl}
{\mathcal F}_\mu(t,v,\dot v,\dots,v^{(\mu+1)})=0,
\end{equation}
which stacks the  equation and its time derivatives
up to level $\mu$ into one large system. We denote partial derivatives of ${\mathcal F}_\mu$ with respect to selected variables~$\zeta$ from $(t,v,\dot v,$ $\ldots, v^{(\mu+1)})$  by ${\mathcal F}_{\mu;\zeta}$,
and the solution set of
the nonlinear algebraic equation associated with the derivative array ${\mathcal F}_{\mu}$
for some integer $\mu$ (considering the variables as well as their derivatives as algebraic variables) by
${\mathbb L}_{\mu}=\{
v_\mu\in{\mathbb I}\times{\mathbb R}^{n+m} \times\ldots
\times{\mathbb R}^{n+m}\mid {\mathcal F}_{\mu}(v_\mu)=0\}$.

The main assumption for the analysis is that the DAE satisfies the following hypothesis, which in the linear case can be proved as a Theorem, see \cite{KunM06}.
\begin{hypothesis}\label{hypns}
Consider the system of nonlinear DAEs (\ref{nlcontrol}). There exist integers $\mu$, $r$, $a$, $d$, and $\nu$
such that ${\mathbb L}_{\mu}$ is not empty and such that for every
$v_\mu^0=(t_0,v_0,\dot v_0,\ldots,v_0^{(\mu+1)})\in{\mathbb L}_{\mu}$
there exists a neighborhood in which
the following properties hold.
\begin{enumerate}
\item
The set ${\mathbb L}_{\mu}\subseteq{\mathbb R}^{(\mu+2)(n+m)+1}$
forms a manifold of dimension $(\mu+2)(n+m)+1-r$.
\item
We have $\rank {\mathcal F}_{\mu;v,\dot v,\ldots,v^{(\mu+1)}}=r$
on ${\mathbb L}_{\mu}$.
\item
We have $\corank {\mathcal F}_{\mu;v,\dot v,\ldots,v^{(\mu+1)}}-
\corank {\mathcal F}_{\mu-1;v,\dot v,\ldots,v^{(\mu)}}=\nu$
on ${\mathbb L}_{\mu}$, where the corank is the dimension of the corange
and the convention is used that  $\corank$ of $ {\mathcal F}_{-1;v}$ is $0$.
\item
We have $\rank {\mathcal F}_{\mu;\dot v,\ldots,v^{(\mu+1)}}=r-a$
on ${\mathbb L}_{\mu}$ such that there exist smooth full rank matrix functions
$Z_2$ and $T_2$
of size $(\mu+1)n \times a$ and $(n+m) \times(n+m-a)$, respectively, satisfying
%
$Z_2^T{\mathcal F}_{\mu;\dot v,\ldots,v^{(\mu+1)}}=0$,
$\rank Z_2^T{\mathcal F}_{\mu;v}=a$, and
$Z_2^T{\mathcal F}_{\mu;v}T_2=0$
%
on ${\mathbb L}_{\mu}$.
\item
We have $\rank {\mathcal F}_{\dot v}T_2=d=n-a-\nu$ on ${\mathbb L}_{\mu}$
such that there exists a smooth full rank matrix function
$Z_1$ of size $n \times d$ satisfying $\rank Z_1^T{\mathcal F}_{\dot v}T_2=d$.
\end{enumerate}
\end{hypothesis}
The smallest $\mu$
for which Hypothesis~\ref{hypns} holds is called the \emph{strangeness-index}
of (\ref{nlcontrol}), see \cite{KunM06}. It generalizes  the concept of \emph{differentiation-index} \cite{BreCP96} to over- and under-determined systems but in contrast to the differentiation-index, ordinary differential equations and purely algebraic equations have $\mu=0$ and for other systems the differentiation-index (if defined) is $\mu+1$, see \cite{KunM06}. The quantity~$\nu$ gives the number of trivial  equations $0=0$ in the system. Of course, these equations can be simply removed and so for our further analysis we assume that $\nu=0$.

If Hypothesis~\ref{hypns} holds
then,   locally (via the implicit function theorem) there exists, see \cite{KunM01,KunM06}, a system  (in the same variables)
%
\begin{eqnarray}
 \hat {\mathcal F}_1(t,v,\dot v)&=&0, \nonumber \\ 
 \hat {\mathcal F}_2(t,v)&=&0,\label{red1}
\end{eqnarray}
%
in which the first $d$ equations $\hat {\mathcal F}_1=Z_1^T{\mathcal F}$ form a (linear)  projection of the original set of equations representing the dynamics of the system, while the second set $\hat {\mathcal F}_2(t,v)=0$  of $a$ equations contains all explicit and hidden constraints and can be used to parameterize the solution manifold and to characterize when an initial condition is consistent. Adding again the  output equation and writing (\ref{red1}) in the original variables we obtain the system
%
\begin{eqnarray}
 \hat {\mathcal F}_1(t,x,\dot x,u)&=&0, \nonumber \\ 
 \hat {\mathcal F}_2(t,x,u)&=&0,\label{red1o}\\
 x(t_0)&=&x^0 \nonumber \\
 y &=& G(t,x,u), \nonumber 
\end{eqnarray}
%
It should be noted that although formally also derivatives of $u$ have been used to form the derivative array, no derivatives of $u$ appear in the regularized system (\ref{red1o}). This has been shown in various contexts \cite{ByeKM97,KunM06,KunMR01} and is  due to the fact that only derivatives of equations where $F_{\mu;u}\equiv 0$  are needed to generate (\ref{red1o}). This means, in particular, that the equations
in  $\hat {\mathcal F}_2(t,x,u)=0$ can be partitioned further into equations that arise from the original system, which include those algebraic equations in the original system which are explicit constraints (in the behavior sense) so that the system can be made to be of differentiation index at most one by a feedback (the part that is impulse controllable or controllable at infinity), and implicit hidden constraints arising from differentiations of equations for which  $F_{\mu,u}\equiv 0$ in the derivative array (the parts that are not impulse controllable).

Using these observations, the regularized system can be (locally in the nonlinear case) written as
%
\begin{eqnarray}
 \hat E_1 \dot x &=&\hat A_1 x +B_1 u, \nonumber \\ 
 0 &=&\hat A_2 x+B_2 u,\label{lred1o}\\
 0 &=&\hat A_3 x, \nonumber \\
 x(t_0)&=&x^0, \nonumber \\
 y &=& C x+Du, \nonumber 
\end{eqnarray}
%
where the third equation that is representing all the hidden algebraic constraints of differentiation index larger that one.
Performing an appropriate (local) change of basis one can identify some (transformed variables) which vanish and the remaining system consisting of the first two equations is of index at most one in the behavior sense, see \cite{ByeKM97,KunM06,KunMR01}
for details. For the first two equations in (\ref{red1o}) and (\ref{lred1o}) one can always find an initial  feedback $u=k(x)+\tilde u$  so that the resulting system is strangeness-free (of differentiation-index one) as a system with input $\tilde u=0$, see \cite{BinMMS15,CamKM11} for a detailed analysis and regularization procedures. In the following we assume that this reinterpretation has been done, so that in the linearization the matrix (functions)
\begin{equation}\label{eq:uode}
\mat{c} \hat E_1 \\ \hat A_2+B_2K \\ \hat A_3 \rix, \
\mat{c} \hat E_1 \\ \hat A_2+B_2K \rix,\
\end{equation}
both  have full row rank, see \cite{KunM06}. Furthermore there exists a (local) partitioning of the variables so that the strangeness-free formulation takes the form
\begin{equation}
\mat{ccc}  \hat E_{11} & \hat E_{12} & \hat E_{13} \\ 0 & 0 & 0 \\ 0 & 0 & 0\rix \mat{c} \dot x_1 \\ \dot x_2 \\ \dot x_3 \rix = \mat{ccc}  \hat A_{11} & \hat A_{12} & \hat A_{13} \\ \hat A_{21} & \hat A_{22} & \hat A_{23}\\  \hat A_{31} & \hat A_{32} & \hat A_{33} \rix \mat{c} x_1 \\ x_2 \\ x_3 \rix
+\mat{c} \hat B_1 \\ \hat B_{2} \\ 0 \rix u
\end{equation}
with the property that $\hat A_{33}$ is invertible and the reduced system obtained by solving for
$x_3$ is strangeness-free (of differentiation index at most one) when setting $u=0$.

The described regularization procedure holds for general DAEs but it does not reflect an available port-Hamiltonian structure. We will now  modify this approach for nonlinear systems with a  pHDAE structure, which (based on the linear time-varying formulation) we define as follows.
\begin{definition}\label{def:local}
Consider a general DAE model  in the form (\ref{nledestate}) and a Hamiltonian
$\mathcal H(x)$ with the property that for a given input $u(t)$ and associated trajectory $x(t)$ the Hessian
$Y_{loc}(t)=\mathcal H_{xx}(x(t))$ can be expressed locally as $E_{loc}(t)^TQ_{loc}(t)$, where
$E_{loc}(t)= F_{\dot x}(t)$,  $F_x(t)=(J_{loc}(t)-R_{loc}(t))Q_{loc}(t) -E_{loc}(t) K_{loc}(t)$, $F_u(t)=B_{loc}(t)-P_{loc}(t)$,
$G_x(t)= (B_{loc}(t)+P_{loc}(t))^TQ_{loc}(t)$, $G_u(t)=S_{loc}(t)+N_{loc}(t)$,
with $E_{loc}(t),A_{loc}(t),Q_{loc}(t),R_{loc}(t)=R_{loc}^T(t), K_{loc}(t) \in C^0({\mathbb I},{\mathbb R}^{n,n})$, $B_{loc}(t),P_{loc}(t)\in C^0({\mathbb I},{\mathbb R}^{n,m})$,  $S_{loc}(t)=S_{loc}^T(t), N_{loc}(t)=-N_{loc}^T(t)\in C^0({\mathbb I},{\mathbb R}^{m,m})$.
 Then the system is called  pHDAE system if the following properties are satisfied:
\begin{itemize}
\item [i)] the (local) differential-algebraic operator
\begin{equation} \label{locskewop}
\mathcal L_{loc}:= Q_{loc}^T(t)E_{loc}(t) \frac{d}{dt} -Q_{loc}^T(t)(J_{loc}(t) Q_{loc}(t) - E_{loc}(t) K_{loc}(t))
\end{equation}
is \emph{skew-adjoint},
\item [ii)] locally  the Hessian $\mathcal H_{xx}$ of the \emph{Hamiltonian} function $\mathcal H(x)$
is bounded from below by a constant;
\item [ii)] locally
%
%
\begin{equation} \label{WdefLoc}
W_{loc}(t)=\left[\begin{array}{lc}
Q_{loc}(t)^T R_{loc}(t) Q_{loc}(t)& Q_{loc}(t)^T P_{loc}(t) \\[1mm]
P_{loc}^T(t) Q_{loc}(t) & S_{loc}(t)
\end{array}\right] \geq 0\ \mbox{\rm for all $t\in \mathbb I$}.
\end{equation}
\end{itemize}
\end{definition}
Clearly standard pH systems of the form (\ref{PHdef}) and linear time-varying systems as in Definition~\ref{def:pHDAE}  directly  fit in this framework. The same holds for multibody systems as in Example~\ref{ex:robot} with the derivative of the constraint
\begin{eqnarray}
     \tilde M(\Theta){\ddot \Theta}+\tilde C(\Theta,{\dot \Theta})+\tilde G(\Theta)&=&\tilde u+\Psi^T\lambda, \nonumber\\
	\Psi(\Theta) \dot \Theta & = & 0. \label{eq:dercon}
\end{eqnarray}
%

\begin{remark}\label{rem:QInon}{\rm There is a lot of choice in the local matrices $Q_{loc}$ and $E_{loc}$ when factoring the Hessian. In some cases we can just choose $Q_{loc}$ to be the identity (see Remark~\ref{remark:pHDAEnoE}), so that $E_{loc}=E_{loc}^T$ defines the Hessian. In other cases one chooses the block-diagonal matrix $E_{loc}=\diag(I_d,0)$ and obtains a semi-explicit formulation of the pHDAE.
However, in general, this freedom should be chosen  to make the system robust to perturbations for simulation and control  methods.
}
\end{remark}
There are multiple reasons why constraints may arise in pH systems. A typical example arises as a limiting situation in a singularly perturbed problem which has pH structure. Typical examples are mechanical multibody systems where small masses are ignored.
\begin{example}
\label{ex:mbseps}{\rm
Finite element modeling of the acoustic field in the interior of a  car, see e.~g. \cite{MehS11}, leads (after several simplifications) to a
large scale constant coefficient differential-algebraic equation system of the form
\[
M\ddot{{p}}+D\dot{{p}}+Kp=B_1u,
\]
where  $p$ is the coefficient vector associated with the pressure in the air and the displacements of the structure,  $B_1u$ is an external force, $M$ is a symmetric positive semidefinite mass matrix, $D$ is a symmetric positive semidefinite matrix, and  $K$ is a symmetric positive definite stiffness matrix. Here $M$ is only semidefinite since
small masses were set to zero, so $M$ is a perturbation of a positive definite matrix. Performing a first order formulation
we obtain  the state equation of a pHDAE system $E\dot z= (J-R)Qz+Bu$, where
\begin{eqnarray*}
E&:=&\mat{cc} M & 0 \\ 0 & I \rix,\;
J:=\mat{cc} 0  &  -I\\ I & 0 \rix,\; R:=
\mat{cc} D  &    \\  & 0 \rix,\; z:=\mat{c}\dot p \\ p \rix,\\
Q&:=&\mat{cc} I & 0 \\ 0 & K \rix,\; B:= \mat{cc}  B_1 \\ 0 \rix,\; P:=0,
\end{eqnarray*}
and the Hamiltonian is
\[
\mathcal H= \frac 12(z^T E^TQz) = \frac 12 (\dot p^T M \dot p +p^T K p).
\]
Note that this model is nonlinear originally, but the simplifications carried out in the modeling process, e.~g. linearization and omitting nonlinear terms with small coefficients leads to a linear model.
}
\end{example}
The other class of examples are systems such as as Example~\ref{ex:robot}, where the dynamics is constrained to a manifold. If the system like in this example has hidden constraints, then the formulation as pHDAE system is not unique because different formulations of the equations and the constraints can be be made. We will come back to this question in Section~\ref{sec:index_reduction}.

As mentioned in the introduction, it is sometimes claimed that port-Hamiltonian DAEs are of differentiation-index at most one (i.~e.,  satisfy Hypothesis~\ref{hypns} with $\mu=0$).  If this would be the case
then in the derivative array with $\mu=0$ the matrix $\mathcal F_{\dot x}$ locally has constant rank $d$ and if $Z_2^T$ is a maximal rank matrix such that (locally) $Z_2^T\mathcal F_{\dot x}=0$ and $Z_1$ is such that it completes $Z_2$ to an invertible matrix $Z=[ Z_1, Z_2]$, then the matrix $ \bar E:=\mat{c} Z_1^T \mathcal F_{\dot x} \\ Z_2^T \mathcal F_{x} \rix$ is locally  invertible.

Let us check this for some of the examples. In Example~\ref{ex:circuit} we have $Z_2^T=\mat{ccc} 0 & 0 & I \rix $ and obtain
\[
\bar E= \mat{ccc} G_c C G_c^T & 0 & 0\\ 0 & L & 0 \\ -G_v^T & 0 & 0 \rix
\]
which is clearly not invertible, except if the last row and column is empty. The same matrix $Z_2$ can be used in Example~\ref{ex:gas}
and yields
\[
\bar E= \begin{bmatrix} M_1 & 0 & 0 \\ 0 & M_2 & 0 \\ 0 & N & 0 \end{bmatrix}
\]
which is also not invertible except if the last row and column is empty. Actually due to the special structure it can be shown that both systems have $\mu=1$, i.~e., differentiation-index two, when the input is chosen to be $0$. The  analysis of Example~\ref{ex:robot} with the original constraint $\Psi(\Theta)$ has $\mu=2$ (differentiation-index three) and the formulation as pHDAE is not straightforward, but using as constraint its derivative  yields $\mu=1$ (differentiation-index two) if $\tilde G \tilde G^T$ is invertible. This replacement corresponds to an index reduction. How to carry out such a regularization for pHDAEs  will be discussed  in Section~\ref{sec:index_reduction}. But let us first (in the next section) analyze in detail the case of differentiation-index one pHDAEs.

\section{PHDAEs of differentiation-index at most one }\label{sec:index_one}
In this section we characterize linear time-varying pHDAE systems of differentiation-index at most one ($\mu=0$).
 In this case Hypothesis~\ref{hypns} implies that the matrix function $E(t)$ has constant rank. Then, see e.~g. Theorem 3.9 in \cite{KunM06}, there exist pointwise orthogonal matrix functions $\tilde U$ and $\tilde V$
such that
\[
\tilde U^TE\tilde V =\mat{cc}E_{11}&0\\0&0\rix=:\tilde E,
\]
where $E_{11}$ is pointwise invertible.  Because $Q^TE$ is real symmetric, setting
\[
\tilde U^TQ\tilde V=\mat{cc}Q_{11}&Q_{12}\\ Q_{21}&Q_{22}\rix,
\]
one has $Q_{11}^T E_{11}=E_{11}^TQ_{11}$ and also $Q_{12}=0$.
Partition in the same way
\begin{eqnarray*}
\tilde U^T J \tilde U &=&\mat{cc}\tilde J_{11} &\tilde J_{12} \\\tilde J_{21} &J_{22} \rix,\quad
\tilde U^T R \tilde U =\mat{cc}\tilde R_{11} &\tilde R_{12} \\\tilde R_{12}^T &R_{22} \rix,\\
\tilde U^T (J -R )\tilde U &=&\mat{cc}\tilde J_{11} &\tilde J_{12} \\\tilde J_{21} &J_{22} \rix
-\mat{cc}\tilde R_{11} &\tilde R_{12} \\\tilde R_{12}^T &R_{22} \rix=:
\mat{cc}\tilde L_{11} &\tilde L_{12} \\\tilde  L_{21} & L_{22} \rix,\\
\tilde
K &=&\tilde V^T(K\tilde V+ \dot{\tilde V}) =\mat{cc}\tilde K_{11} &K_{12} \\\tilde K_{21} &K_{22} \rix,
\end{eqnarray*}
%
%
%
Since the system has differentiation-index at most one, the block $L_{22}Q_{22}$ either does not occur (in this case we have an implicitly defined standard pH  system) or it must be pointwise invertible, see \cite{KunM06}, i.~e., both
$L_{22}$ and $Q_{22}$ are pointwise invertible. Let  $U=\tilde UT$, where
\[
T:=\mat{cc}I&0\\T_{21}&I\rix,\qquad T_{21}=-L_{22}^{-T} (\tilde L_{12}-E_{11}K_{12}Q_{22}^{-1})^T.
\]
Then a transformation of the original pHDAE with $U$  and $\tilde V$ yields a transformed pHDAE system,
where $\tilde K$ is defined above,
%
\begin{eqnarray*}
\tilde E&=&U^T EV=\tilde U^T E\tilde V, \quad \tilde Q=U^{-1}Q\tilde V=\mat{cc}Q_{11}&0\\\tilde Q_{21}&Q_{22}\rix,\quad
\tilde S=S, \quad \tilde N=N,\\
\tilde J&=&U^TJU=\mat{cc}J_{11}& J_{12}\\ J_{21}&J_{22}\rix,\quad
\tilde R=U^TRU=\mat{cc}R_{11}& R_{12}\\ R_{12}^T&R_{22}\rix,\\
\tilde L&=&\tilde J-\tilde R=\mat{cc}J_{11}& J_{12}\\ J_{21}&J_{22}\rix-
\mat{cc}R_{11}& R_{12}\\ R_{12}^T&R_{22}\rix=\mat{cc}L_{11}& L_{12}\\ L_{21}&L_{22}\rix,
\end{eqnarray*}
 and
\[
\tilde L\tilde Q-\tilde E\tilde K
=\mat{cc}L_{11}Q_{11}+L_{12}\tilde Q_{21}-E_{11}\tilde K_{11}&0\\
L_{21}Q_{11}+L_{22}\tilde Q_{21}&L_{22}Q_{22}\rix,
\]
 %
%
%
%
i. e.,
\begin{equation}
\label{eq:1HZ}
  ({J}_{12} - {R}_{12}) Q_{22} - E_{11} K_{12} =0.
\end{equation}
Performing another change of basis to make $\tilde{Q}$ (block) diagonal with a transformation matrix
\[
   \tilde T:=  \left[ \begin{array}{cc} I & 0\\ - Q_{22}^{-1} \tilde{Q}_{21} & I \end{array} \right],
\]
then setting $V= \tilde V\tilde T$ and transforming the original system with $U,V$ we obtain that any
pHDAE of {{ differentiation-}}index one can be transformed to the form
\begin{eqnarray}
\nonumber
\mat{cc}E_{11}&0\\0&0\rix \mat{c}\dot x_1\\\dot{x}_2 \rix
&=&\left(\mat{cc} L_{11}&L_{12}\\ L_{21}&L_{22}\rix\mat{cc}Q_{11}&0\\ 0 &Q_{22}\rix
-\mat{cc}E_{11}K_{11}&E_{11}K_{12}\\0&0\rix\right)\mat{c}x_1\\x_2\rix\\
\label{eq:2HZ}
&& \qquad+ \mat{c}B_1-P_1\\B_2-P_2\rix u, \\
\nonumber
y&=& \mat{cc}(B_1+P_1)^T&(B_2+P_2)^T\rix\mat{cc}Q_{11}&0\\ 0 &Q_{22}\rix\mat{c}
x_1\\x_2\rix + (S+N) u,
\end{eqnarray}
where (\ref{eq:1HZ}) holds, with
\[
K_{11}=\tilde K_{11}-K_{12}Q_{22}^{-1}\tilde Q_{21}^{-1},\qquad
\mat{cc}B_1&P_1\\B_2&P_2\rix =U^T\mat{cc}B&P\rix.
\]
%
%
%
Following Theorem~\ref{pHDAEinv} this transformation will not change the Hamiltonian, and
(\ref{eq:2HZ}) is still a pHDAE of index at most one.
Note that these transformations should not be performed in a numerical integration or control design technique, since the inversion of the matrices $Q_{22}$ and $L_{22}$ may be highly ill-conditioned. However, from an analytic point of view have the following theorem.
%
\begin{theorem}\label{redpHDAEthm}
Suppose that the pHDAE system (\ref{pHDAE}) is of differentiation-index at most one (i.~e.\ it satisfies Hypothesis~\ref{hypns} with $\mu=0$) that $\nu=0$,  and that $E(t)$ has constant rank.  Assume further that the system is transformed to the form (\ref{eq:1HZ})--(\ref{eq:2HZ}). Then for any input function $u$ and  $x_1(t_0)=x_{1,0}$ the first component of the solution and the output of (\ref{eq:2HZ}) are given by reduced implicit pHDAE system
\begin{eqnarray}
\nonumber
E_{11} \dot x_1&=&[(J_{11}- R_{11}) Q_{11} - E_{11} K_{11} ]x_1+(\hat B-\hat P)u, \qquad x_1(t_0)=x_{1,0},\\
\label{redpHDAE}
y &=&(\hat B+\hat P)^T Q_{11} x_1+(\hat S+\hat N)u,
\end{eqnarray}
with Hamiltonian $\hat{\mathcal H}(x_1)=\frac{1}{2}x_1^T  Q_{11}^T E_{11} x_1
={\mathcal H}(x)$, and coefficients
\begin{eqnarray*}
\hat B&=& B_1-\frac{1}{2}( J_{21}^T- R_{12})L_{22}^{-T}(B_2+P_2),\\
\hat P &=&  P_1-\frac{1}{2}( J_{21}^T- R_{12})L_{22}^{-T}(B_2+P_2),\\
\hat S&=&S-\frac{1}{2}[(B_2+P_2)^TL_{22}^{-1}(B_2-P_2)+(B_2-P_2)^{T}L_{22}^{-T}(B_2+P_2)],\\
\hat N&=&N-\frac{1}{2}[(B_2+P_2)^TL_{22}^{-1}(B_2-P_2)-(B_2-P_2)^{T}L_{22}^{-T}(B_2+P_2)].
\end{eqnarray*}
Furthermore, the second part of the state $x_2$ is uniquely determined by the algebraic constraint
\begin{equation}\label{expx2}
L_{22}Q_{22}x_{2}=-L_{21}Q_{11}x_1 -(B_2-P_2)u,
\end{equation}
and there is a consistency constraint for the initial condition
\begin{equation}
\label{eq:7HZ}
  L_{22}(t_0)Q_{22}(t_0)x_{2}(t_0)=-[(L_{21}(t_0))Q_{11}(t_0)]x_1(t_0)
-(B_2(t_0)-P_2(t_0))u(t_0).
\end{equation}
\end{theorem}
\proof
Equation (\ref{expx2}) follows directly from the second state equation in (\ref{eq:2HZ}).
Since $\hat{B}-\hat{P}= B_1-P_1$, we see that $x_1$ satisfies the state equation in (\ref{redpHDAE}).
The output equation is obtained directly by substituting (\ref{expx2}) in the output equation of (\ref{eq:2HZ}).

It remains to prove that (\ref{redpHDAE}) is port-Hamiltonian. Since (\ref{eq:2HZ}) is a pHDAE system, it follows that
\eq{eq:5HZ}
Q_{11}^T E_{11}= E_{11}^T Q_{11}
\en
and
%
%
\begin{eqnarray}
\nonumber
\frac{d}{dt}Q_{11}^TE_{11}
&=&Q_{11}^T[ E_{11}K_{11} - J_{11}Q_{11}]
+[ E_{11}K_{11} - J_{11}Q_{11}]^T Q_{11},\\
\label{eq:4HZ}
0&=&-Q_{11}^T( J_{12} +\ J_{21}^T)Q_{22} +Q_{11}^TE_{11}K_{12} ,\\
\nonumber
0&=& Q_{22}^T J_{22} Q_{22} + Q_{22}^T J_{22}^T Q_{22}.
\end{eqnarray}
%
Combining (\ref{eq:5HZ}) with the first equation of (\ref{eq:4HZ}) gives that the operator $Q_{11}^T E_{11}\frac{d}{dt} - Q_{11}^T( J_{11} Q_{11} - E_{11} K_{11})$ is skew-adjoint.

Furthermore, since $\hat{S}$ is symmetric and $\hat{N}$ is skew-symmetric,  system (\ref{redpHDAE}) is of the form (\ref{pHDAE}), and thus Theorem \ref{thm:dissin} gives that (\ref{eq:bal}) is satisfied. So
\begin{equation}
\label{eq:8HZ}
\frac{d}{dt}\hat{\mathcal H}(x_1)=
\frac{d}{dt} x_1^T Q_{11}^T E_{11} x_1 =y^Tu-\mat{c}x_1\\ u\rix^T \hat{W} \mat{c}x_1\\u\rix
\end{equation}
with
\[
  \hat W=\mat{cc}Q_{11}^T R_{11} Q_{11}&Q_{11}^T\hat P\\
   \hat P^TQ_{11}&\hat S\rix.
\]
On the other hand, since (\ref{pHDAE}) is a pHDAE system, we have that
\begin{equation}
\label{eq:9HZ}
  \frac{d}{dt}{\mathcal H}(x)= y^Tu - \mat{c} x\\ u \rix^T W \mat{c} x\\ u \rix,
\end{equation}
where
\[
  W=\mat{ccc}Q_{11}&0&0\\ 0 &Q_{22}&0\\0&0&I\rix^T
  \mat{ccc}  R_{11}& R_{12}& P_1\\ R_{12}^T&R_{22}&P_2\\ P_1^T&P_2^T&S\rix
\mat{ccc}Q_{11}&0&0\\ 0 &Q_{22}&0\\0&0&I\rix.
\]
We know that for the same input and initial state with $x_2(t_0)$ satisfying (\ref{eq:7HZ}) the solutions of the two systems are the same, and furthermore we have that
\[
  {\mathcal H}(x) = x^T Q^T E x = x_1^T Q_{11}^T E_{11} x_1 = \hat{\mathcal H}(x_1),
\]
and that (\ref{expx2}) holds. Thus, from (\ref{eq:8HZ}) and (\ref{eq:9HZ}) we obtain that
\begin{equation}
\label{eq:10HZ}
  \mat{c}x_1\\u\rix^T \hat{W} \mat{c}x_1\\u\rix = \mat{c} x_1\\ x_2\\ u \rix^T W \mat{c} x_1\\ x_2 \\u \rix=  \mat{c}x_1\\u\rix^T W_X \mat{c}x_1\\u\rix,
\end{equation}
where $W_X= X^T  WX $ with
\[
X=\mat{cc}I&0\\-Q_{22}^{-1}(L_{22}^{-1}( J_{21}-R_{12}^T)Q_{11})&-Q_{22}^{-1}L_{22}^{-1}(B_2-P_2)\\0&I\rix.
\]
Since (\ref{eq:10HZ}) has to hold for all $x_1$ and $u$, we find that $\hat{W} = W_X$, which could also be obtained by straightforward (but tedious) calculation.
Since $W$ is symmetric positive semidefinite, so is $\hat W$, and hence the reduced system in  $x_1$ is still port-Hamiltonian with Hamiltonian $\hat{\mathcal H}(x_1)$.
\eproof
Note that for the numerical integration or in the control context, as for general DAEs, it is sufficient to carry out the transformation with pointwise orthogonal $\tilde U$ from the left and the insertion of $I=\tilde U \tilde U^T$ before $Q$. In this way a differentiation of a computed transformation matrix can be avoided and the pHDAE structure is preserved nonetheless. However, no explicit separation of the parts $x_1$ and $x_2$ would be obtained in this way and this separation has to be carried out by the numerical solver in the context of the numerical integration method.

\begin{remark}\label{rem:nl}{\rm
For nonlinear pHDAE systems with differentiation index at most one ($\mu=0$), the corresponding local result follows directly via the implicit function theorem
and application of Theorem~\ref{redpHDAEthm} to the linearization as in Definition~\ref{def:local}.}
\end{remark}

\section{Regularization of higher index pHDAE systems}\label{sec:index_reduction}
In this section we discuss  how to  modify the regularization procedure discussed for general DAEs in Section~\ref{sec:nonl} to preserve the pHDAE structure.
Let us first consider the linear time-varying case (\ref{pHDAE}) and set $L=J-R$. Suppose that
the state equation with $u=0$ already satisfies Hypothesis~\ref{hypns}, i.~e.,  as discussed in Section~\ref{sec:nonl}, no  reinterpretation of variables or initial feedbacks are necessary. It has been shown in \cite{ByeKM97} that the extra constraint equations (hidden constraints) that arise from derivatives are uncontrollable, because otherwise the index reduction could have been done via feedback. This means that these extra constraint equations are of the form
$\hat A_3 x=0$ which corresponds to $\hat F_3(t,x)=0$ in the nonlinear case. We add just these constraint equations to our original
pHDAE and obtain an  overdetermined strangeness-free system. Note again that under our assumptions the explicit algebraic constraints are included in the first two equations in (\ref{red1o}), resp. (\ref{lred1o}).

Let us make the weak assumption that $E(t)$ has constant rank. This is a restriction that however holds in all examples that we have encountered so
far, and it can be removed by considering the system in a piecewise fashion, see \cite{KunM06}. Then there exist real orthogonal matrix functions $U, V_1\in C^1(\mathbb I, \mathbb R^{n,n})$ such that
\[
U_1^T E V_1 =\mat{cc} \tilde E_{11} & 0\\
0 & 0 \rix
\]
with pointwise invertible $\tilde E_{11}$.
Perform a transformation of the pHDAE (\ref{pHDAE}) as in  Theorem~\ref{pHDAEinv} and also form $\hat A_3 V= \mat{cc} \hat A_{31} & \hat A_{32} \rix$ partitioned accordingly. By the property that $\hat A_3$ contains all the high index constraints it follows that
$\hat A_{23}$ has full row rank for all $t\in \mathbb I$, and hence there exists a real orthogonal matrix function $V_2$ such that $\hat A_{32} V_2= \mat{cc} 0 & A_{33} \rix$ with $A_{33}$ pointwise invertible. Performing a change of variables of the pHDAE with
\[
V:=V_1 \mat{cc} I & 0 \\ 0 & V_2 \rix \mat{ccc} I & 0 & 0 \\ 0 & I & 0 \\ -A_{33}^{-1}\hat{A}_{31}  & 0 & I \rix
\]
we obtain  a pHDAE of the form
\begin{eqnarray}
\nonumber
\lefteqn{\left[\begin{array}{ccc} \tilde E_{11} & 0 &0\\
0 & 0 & 0  \\ 0 & 0 & 0\end{array}\right] \left[\begin{array}{c} \dot x_1 \\ \dot x_2 \\ \dot x_3\end{array}
\right]
= \tilde L\mat{ccc} \tilde Q_{11} &\tilde Q_{12} & \tilde Q_{13} \\\tilde Q_{21} &\tilde Q_{22} &\tilde Q_{23}   \\  \tilde Q_{31} & \tilde Q_{32} & \tilde Q_{33}  \rix \mat{c} x_1 \\ x_2 \\ x_3\rix} \quad \\
 &&\qquad - \mat{ccc}\tilde  E_{11} & 0 &0\\
0& 0 & 0  \\ 0& 0 & 0\rix \mat{ccc}  K_{11}&  K_{12} &    K_{13} \\   K_{21}&  K_{22}  &   K_{23} \\  K_{31}&  K_{32} &    K_{33} \rix
 \mat{c} x_1 \\ x_2 \\ x_3\rix \nonumber \\
&& \qquad+
\mat{c} \tilde B_1 - \tilde P_1 \\   \tilde B_2 -  \tilde P_2 \\ \tilde B_3 -\tilde P_3 \rix u,\label{eq:38}\\
y&=& \mat{ccc} (\tilde B_1+\tilde P_1)^T & (\tilde B_2+\tilde P_2)^T & (\tilde B_3+\tilde P_3)^T \rix \mat{ccc}\tilde Q_{11} &\tilde Q_{12} & \tilde Q_{13} \\\tilde Q_{21} &\tilde Q_{22} &\tilde Q_{23}   \\  \tilde Q_{31} & \tilde Q_{32} & \tilde Q_{33}  \rix \mat{c} x_1 \\ x_2 \\ x_3\rix\nonumber \\
 && + (S+N) u,
\nonumber
\end{eqnarray}
where $\tilde K=(V^TKV+\dot V) $, $\tilde L= LV$, together with the constraint $0=A_{33}x_3$, i.~e. $x_3=0$.

Assuming further that the matrix function
\[
\mat{cc} \tilde Q_{11} &\tilde Q_{12} \\\tilde Q_{21} &\tilde Q_{22}    \\  \tilde Q_{31} & \tilde Q_{32}  \rix
\]
has constant  rank, there exists a pointwise real  orthogonal matrix function $U_2$ such that
\[
U_2^T \mat{ccc} \tilde Q_{11} &\tilde Q_{12} & \tilde Q_{13} \\\tilde Q_{21} &\tilde Q_{22} &\tilde Q_{23}   \\  \tilde Q_{31} & \tilde Q_{32} & \tilde Q_{33}  \rix =\mat{ccc} Q_{11} & Q_{12} &  Q_{13} \\Q_{21} & Q_{22} & Q_{23}   \\  0 & 0& Q_{33}  \rix
\]
Transforming  the pHDAE (\ref{eq:38}) with $U_2^T$ we get a  pHDAE of the form
\begin{eqnarray}
\nonumber
\lefteqn{\mat{ccc} E_{11} & 0 &0\\
E_{21} & 0 & 0  \\ E_{31} & 0 & 0\rix \mat{c} \dot x_1 \\ \dot x_2 \\ \dot x_3\rix
= \mat{ccc}L_{11}&L_{12} &  L_{13} \\ L_{21}&L_{22}  & L_{23}\\ L_{31}&L_{32} &  L_{33} \rix \mat{ccc} Q_{11} & Q_{12} &  Q_{13} \\Q_{21} & Q_{22} & Q_{23}   \\  0 & 0& Q_{33}  \rix \mat{c} x_1 \\ x_2 \\ x_3\rix}\qquad \\
 &&\qquad -\mat{ccc}\tilde  E_{11} & 0 &0\\
 \tilde E_{21}& 0 & 0  \\ \tilde E_{31}& 0 & 0\rix \mat{ccc}  K_{11}&  K_{12} &    K_{13} \\   K_{21}&  K_{22}  &   K_{23} \\  K_{31}&  K_{32} &    K_{33} \rix
 \mat{c} x_1 \\ x_2 \\ x_3\rix \nonumber \\
&& \qquad+
\mat{c} \tilde B_1 - \tilde P_1 \\   \tilde B_2 -  \tilde P_2 \\ \tilde B_3 -\tilde P_3 \rix u,\label{eq:39}\\
y&=& \mat{ccc} (\tilde B_1+\tilde P_1)^T & (\tilde B_2+\tilde P_2)^T & (\tilde B_3+\tilde P_3)^T \rix\mat{ccc}\tilde Q_{11} &\tilde Q_{12} & \tilde Q_{13} \\\tilde Q_{21} &\tilde Q_{22} &\tilde Q_{23}   \\  \tilde Q_{31} & \tilde Q_{32} & \tilde Q_{33}  \rix \mat{c} x_1 \\ x_2 \\ x_3\rix\nonumber \\
 && + (S+N) u,
\nonumber
\end{eqnarray}
together with the constraint $0=x_3$.

By Theorem~\ref{pHDAEinv}, system (\ref{eq:39}) is still a pHDAE system, and the solution of the overdetermined system (\ref{eq:39}) together with $x_3=0$ is the same as that of (\ref{eq:38}) and the Hamiltonian is unchanged.  Since the resulting system is
still port-Hamiltonian, using that $x_3=0$, we have that the subsystem given by the first two block rows (\ref{eq:39}) together with output equation is an index at most one phDAE which has the form
\begin{eqnarray}
\nonumber
\mat{cc} E_{11} & 0 \\
E_{21} & 0 \rix \mat{c} \dot x_1 \\ \dot x_2 \rix
&=& \mat{cc}L_{11}&L_{12}  \\ L_{21}&L_{22} \rix \mat{cc} Q_{11} & Q_{12} \\Q_{21} & Q_{22}   \rix \mat{c} x_1 \\ x_2 \rix\\
 &-&\mat{cc} E_{11} & 0\\
E_{21} & 0 \rix \mat{cc}  K_{11}&  K_{12}  \\   K_{21}&  K_{22}  \rix
 \mat{c} x_1 \\ x_2 \rix +
\mat{c} \tilde B_1 - \tilde P_1 \\   \tilde B_2 -  \tilde P_2  \rix u,\label{eq:40}\\
y&=& \mat{cc} (\tilde B_1+\tilde P_1)^T & (\tilde B_2+\tilde P_2)^T  \rix\mat{cc}\tilde Q_{11} &\tilde Q_{12}  \\\tilde Q_{21} &\tilde Q_{22}  \rix \mat{c} x_1 \\ x_2 \rix\nonumber \\
 &+& (S+N) u,
\nonumber
\end{eqnarray}
To this system we can apply the results of the previous section and obtain that the system can be further reduced to an implicit standard pH system.
\begin{example}{\rm
Consider again the semidiscretized Example~\ref{ex:gas}. It has been shown in \cite{EggKLMM17} that for a (permuted) singular value decomposition (SVD) of $N^T$
\[
N^T =U_N^T \begin{bmatrix} 0  \\  \Sigma \end{bmatrix} V_N,
\]
with real orthogonal matrices $U_N,V_N$ and a nonsingular diagonal matrix $\Sigma \in \mathbb R^{n_3,n_3}$.
Scaling the second row of \eqref{eq:ph1} with $U_N$ and setting $x_2=V_N \begin{bmatrix} x_{2,2}^T & x_{2,3}^T \end{bmatrix}^T$, as well as $x_{2}^0=V_N \begin{bmatrix} {x_{2,2}^0}^T & {x_{2,3}^0}^T \end{bmatrix}^T$ we obtain a transformed system
\begin{equation}\label{eq:splitsystem}
\begin{bmatrix} M_1 & 0 & 0 &0\\ 0 & M_{2,2} &M_{2,3} & 0 \\ 0 & M_{2,3}^T & M_{3,3} & 0 \\ 0 &0 & 0 & 0 \end{bmatrix}
\frac{d}{dt} \begin{bmatrix} x_1 \\ x_{2,2}\\ x_{2,3} \\ x_3 \end{bmatrix}
+\begin{bmatrix} 0 & G_{1,2} & G_{1,3} &0 \\ -G_{1,2}^T & D_{2,2} & D_{2,3} & 0\\ -G_{1,3}^T & D_{2,3}^T & D_{3,3} & -\Sigma \\ 0 & 0 &\Sigma  & 0\end{bmatrix}\begin{bmatrix} x_1 \\ x_{2,2}\\ x_{2,3} \\ x_3 \end{bmatrix}= \begin{bmatrix} 0 \\ B_{2,2} \\ B_{3,2} \\ 0 \end{bmatrix} u.
\end{equation}
It follows immediately that $x_{2,3}=0$, which is the uncontrollable (index two) constraint in the DAE that in particular the initial condition $x^{0}_{2,3}$ has to satisfy.
The vectors $x_1,x_{2,2}$ are solutions of the implicit ordinary pH system
\begin{equation}\label{eq:redodesystem}
\begin{bmatrix} M_1 & 0\\  0 & M_{2,2}  \end{bmatrix}
\frac{d}{dt} \begin{bmatrix} x_1 \\ x_{2,2}  \end{bmatrix}
+\begin{bmatrix} 0 &  G_{1,2}  \\ -G_{1,2}^T & D_{2,2} \end{bmatrix}\begin{bmatrix} x_1 \\ x_{2,2}  \end{bmatrix}= \begin{bmatrix} 0 \\  B_{2,2} \end{bmatrix} u,
\end{equation}
with initial conditions $x_1(0)=x^{0}_1$,  $x_{2,2}(0)=x_{2,2}^0$, so they are well-defined continuously differentiable functions for any piecewise continuous $u$ and any choice of the initial conditions.

Finally we get the component $x_3$ (the Lagrange multiplier) via
\begin{equation}\label{eq:x21}
x_3=\Sigma^{-1}(M_{2,3}^T \frac{d}{dt} x_{2,2}-G_{1,3}^T x_1 +D_{2,3}^T x_{2,2}- B_{3,2} u),
\end{equation}
and this is the implicit index one constraint in the DAE. Since both type of (the explicit and the hidden) constraints have to be satisfied for the initial condition, it means that the transformed initial condition also has to satisfy the consistency condition
\begin{equation}\label{eq:consis}
x_3(0)=\Sigma^{-1}(M_{2,3}^T \frac{d}{dt} x_{2,2}(0)-G_{1,3}^T x_1(0) +D_{2,3}^T x_{2,2}(0)- B_{3,2} u(0))
\end{equation}
Condition~\eqref{eq:consis} leads to a relationship between the input $u$ and the state
at $t=0$, which is a constraint that has to be satisfied to have a classical solution.
Furthermore, we see immediately that to obtain a continuous $x_3$ the function $B_{3,2} u$ has to be continuous and $u$ has to be such that $B_{3,2} u$ leads to a continuous $M_{2,3}^T \frac{d}{dt} x_{2,2}$. The implicit ordinary pH system~\eqref{eq:redodesystem} describes the dynamics of the system, while the other two equations describe the constraints.
}
\end{example}

\begin{remark}\label{rem:nlhi}{\rm
For nonlinear pHDAE systems satisfying Hypothesis~\ref{hypns} with $\mu>0$, the corresponding local result follows directly via linearization and the implicit function theorem.}
\end{remark}
\section*{Conclusion}
A new definition of port-Hamiltonian descriptor systems has been derived. It has been shown that this formulation is valid also for DAEs of differentiation-index larger than one, and it has been demonstrated that under some (local) constant rank assumption any such pHDAE can be reformulated as an implicitly defined standard PH system plus an algebraic constraint that
describes the manifold where the dynamics of the system takes place and that also describes the consistent initial conditions.
As for standard DAEs the reformulated system is well suited for numerical integration and control, since all constraints are available.

\section*{Acknowledgments}
We acknowledge many interesting discussions with Robert Altmann and Philipp Schulze from TU Berlin
and Arjan Van der Schaft from RU Groningen. The first author has been supported by {\it Einstein Foundation Berlin},
through an Einstein Visiting Fellowship.
The second author has been supported by Deutsche Forschungsgemeinschaft for Research support via Project A02 in CRC 1029 TurbIn and by {\it Einstein Foundation Berlin} within the Einstein Center ECMath.

\end{document}